%
%
%  AUTHORS:  FRANC FORSTNERIC AND JASNA PREZELJ
%
%  TITLE:    EXTENDING HOLOMORPHIC SECTIONS FROM COMPLEX SUBVARIETIES
%
%  Final version: July 4, 2000
%
%  Mathematische Zeitschrift 
%
%___________________________________
%
%   This is the Macro file.
%
%
\scrollmode \magnification=\magstep1
\parskip=\smallskipamount

\def\demo#1:{\par\medskip\noindent\it{#1}. \rm}
\def\ni{\noindent}               % noindent
     % new indented par in horizontal mode
 % new unindented par in horizontal mode
\def\ll{\leftline}
\def\cl{\centerline}

\def\begin{\ll{}\vskip 10mm\nopagenumbers}  % beginning of the paper
\def\pn{\footline={\hss\tenrm\folio\hss}}   % pagenumbers at bottom
\def\ii#1{\itemitem{#1}}

%
%   \beginsection
%
\outer\def\beginsection#1\par{\bigskip
  \message{#1}\leftline{\bf\&#1}
  \nobreak\smallskip\vskip-\parskip\noindent}

%
%  \proclaim
%
\outer\def\proclaim#1:#2\par{\medbreak\vskip-\parskip
    \noindent{\bf#1.\enspace}{\sl#2}
  \ifdim\lastskip<\medskipamount \removelastskip\penalty55\medskip\fi}

\def\endpr{\hfill $\spadesuit$ \medskip}

%
%
%  Roman capital characters
%
%
               % roman D in math mode

               % roman T
               % roman N in math

%
%
%  Boldface capital characters
%
%

\def\C{{\bf C}}

\def\R{{\bf R}}

\def\Z{{\bf Z}}

% \def\O{{\bf O}}

%
%
% Caligraphic capital characters
%
%
\def\cA{{\cal A}}

\def\cC{{\cal C}}

\def\cH{{\cal H}}

\def\cJ{{\cal J}}
\def\cK{{\cal K}}

\def\cO{{\cal O}}

\def\cU{{\cal U}}

%
%
%  Small Greek letters in Math mode
%
\def\a{\alpha}
\def\b{\beta}
\def\g{\gamma}
\def\d{\delta}
\def\e{\epsilon}

\def\l{\lambda}

\def\c{\chi}
%\def\o{\omega}

%\def\O{\Omega}

%
%
%  Miscellaneous symbols
%
%
\def\bar{\overline}              % conjugate
\def\bs{\backslash}              % backslash
                % partial derivative
\def\dibar{\bar\partial}         % di-bar derivative
\def\hra{\hookrightarrow}

%
%  Symbols
%
             % disc
        % closed disc
               % boundary of the disc
\def\cn{{\bf C}^n}

\def\c*{{\bf C}^*}
\def\Aut{{\rm Aut}}                         % roman Aut in math mode
              % sypmlectic automorphisms

\def\wt{\widetilde}

%
% Abbreviations
%
                    % dimension
\def\holo{holomorphic}                   % holomorphic
                  % automorphism
                  % homomorphism
               % analytic subset
                % homeomorphism
                    % continuous
\def\nbd{neighborhood}                   % neighborhood
\def\psc{pseudoconvex}                   % pseudoconvex
\def\spsc{strongly\ pseudoconvex}        % strongly psc
                   % real-analytic
\def\psh{plurisubharmonic}               % plurisubharmonic

                    % totally real
             % polynomially convex
\def\hc#1{${\cal H}(#1)$-convex}         % holo convex
           % holomorphic function
\def\ss{\subset\!\subset}                % relatively compact subset

\def\supp{{\rm supp}\,}                  % support
              % C^n equivalent

\def\hvf{holomorphic vector field}
\def\hvb{holomorphic vector bundle}

%_________________________________________________________

\begin
\cl{\bf EXTENDING HOLOMORPHIC SECTIONS}
\cl{\bf FROM COMPLEX SUBVARIETIES}
\bigskip
\cl{Franc Forstneri\v c and Jasna Prezelj}
\bigskip\medskip\rm

%
%
%  INTRODUCTION
%
%
\beginsection 1. Introduction and results.

In this paper we treat the classical problem of extending \holo\ 
mappings and sections from closed complex subvarieties of a complex
manifold. Our main results (theorems 1.1 and 1.4) extend 
those of Grauert [Gr2] and Cartan [Car]; results of this type 
are commonly referred to as the `Oka-Grauert principle'
on Stein manifolds.  Our methods are similar to those developed 
in the papers of Gromov [Gro], Henkin and Leiterer [HL2] and 
the authors [FP1], [FP2]. The main addition here is the interpolation 
of a given holomorphic section on a complex subvariety of a Stein manifold. 
Our globalization scheme follows very closely the one developed in [FP2].

To state our first result we recall the notion of a (dominating)
spray introduced by M.\ Gromov ([Gro], sec.\ 0.5). Given a
\hvb\ $p\colon E\to Y$ over a complex manifold $Y$, 
we denote by $0_y\in E_y=p^{-1}(y)$	 the zero element in the 
fiber $E_y$ and we observe that $E_y$ is a $\C$-linear subspace 
of the tangent space $T_{0_y} E$.

%
%  Sprays
%
\proclaim Definition 1:
A {\bf spray} on a complex manifold $Y$ is a \hvb\ $p\colon E\to Y$,
together with a holomorphic map $s\colon E\to Y$, such that
for each $y\in Y$, $s(0_y)=y$ and the derivative $ds\colon
T_{0_y}E \to T_y Y$ maps $E_y$ surjectively onto $T_y Y$.

%
%  Main theorem for mappings
%
\proclaim 1.1 Theorem: Let $X$ be a Stein manifold and $Y$ a
complex manifold which admits a spray. Then for every closed
complex subspace $X_0\subset X$ and every continuous map
$f_0\colon X\to Y$ whose restriction to $X_0$ is \holo\ on $X_0$
there exists a \holo\ map $f_1\colon X\to Y$ such that
$f_1|_{X_0}= f_0|_{X_0}$. Moreover, $f_1$ can be obtained
from $f_0$ by a homotopy $f_t\colon X\to Y$ ($t\in [0,1]$)
which is fixed on $X_0$.

For maps of Stein manifolds into complex Lie groups or
complex homogeneous spaces this was proved by Grauert
([Gra2], [Gra3]) and Cartan ([Car], Th\'eor\`eme 1 bis).
The validity of theorem 1.1 (and of theorem 1.4 below)
was asserted by Gromov ([Gro], sec.\ 2.9.C, p.\ 877),
but very few details were provided there.
The corresponding results with $X_0=\emptyset$
can be found in [Gro], and complete proofs are given in [FP1]
and [FP2]. (For the general theory of Stein manifolds
and Stein spaces we refer to [GuR] or [GRe].)

The best source of examples of complex manifolds with sprays
is the following. Let $V_1,\ldots, V_q$ be $\C$-complete \hvf s
on a complex manifold $Y$ which span the tangent space $T_y Y$
at each point $y\in Y$. Denote by $\theta^j_t$ the flow of $V_j$.
Then the map $s\colon Y\times \C^q\to Y$,
$s(y,t_1,\ldots,t_q)= \theta_{t_1}^1\circ\cdots\circ\theta_{t_q}^q(y)$
for $y\in Y$ and $(t_1,\ldots,t_q) \in \C^q$, is a spray on $Y$
(see [Gro] and [FP1]). Such sprays exist on complex Lie groups
and homogeneous spaces (take left invariant \hvf s spanning
the Lie algebra), and on spaces $Y=\C^n\bs \Sigma$ where
$\Sigma$ is an affine algebraic subvariety of codimension
at least two (see [FP2]).

Theorem 1.1 reduces the \holo\ extension problem to a
topological one. If the initial map $f_0$ is
\holo\ in an open \nbd\ of $X_0$ in $X$, the homotopy
in theorem 1.1 can be chosen such that all sections $f_t$ are
\holo\ in a \nbd\ of $X_0$ (independent of $t$) and they agree
with $f_0$ to a prescribed finite order along $X_0$
(see theorem 1.4). A local \holo\ extension of $f_0$
exists without any condition on $X$ and $Y$ provided that
the subspace $X_0$ is Stein (proposition 1.3).
Theorem 1.1 is a special case of theorem 1.4 below
which gives extension results for sections of \holo\
submersions onto Stein manifolds.

We have already pointed out that theorem 1.1 applies
to maps of Stein manifolds into $\C^n\bs \Sigma$ where
$\Sigma$ is an affine algebraic subvariety of codimension
at least two. This does not hold for general analytic subvarieties,
independently of their codimension.

\pn

\demo Example 1: For each $n\ge 1$ there exists a discrete
set $\Sigma\subset \C^n$ such that theorem 1.1 fails for maps of
Stein manifolds into $Y=\C^n\bs\Sigma$. (Note however that
theorem 1.1 does hold when $\Sigma$ is a discrete set in $\C^n$,
$n\ge 2$, which is tame in the sense of Rosay and Rudin [RRu],
since the complement of such a set admits a spray.)

We obtain such examples as in the proof of theorem 1.6 (b) in [FP2].
Let $n=2$ for simplicity. By [RRu] there is a discrete set
$\Sigma\subset\C^2$ such that the space $Y=\C^2\bs\Sigma$ is
 volume hyperbolic, in the sense that any entire holomorphic
map $g\colon \C^n \to Y$ has complex rank at most one at each point.
The same is then true for holomorphic maps $f\colon X\to Y$ from
any Stein manifold $X$ whose universal cover is a complex Euclidean space.
We showed in [FP2] that consequently any such map is homotopic
to a constant map $X\to {\rm point}$.

The set $X=(\C^*)^2 =\C^2\bs\{zw=0\}$ is a Stein manifold which
is covered by $\C^2$. Choose a complex line
$\Lambda\subset \C^2$ that intersects the set
$\{zw=0\} \subset \C^2$
in precisely two points, say $p=(1,0)$ and $q=(0,1)$, and
take $X_0=\Lambda\cap (\C^*)^2 =\Lambda \bs\{p,q\}$; this twice
punctured complex line is a closed complex submanifold of $X$. We
may choose the discrete set $\Sigma \subset \C^2$ as above such
that $\Sigma\cap\Lambda =\{p,q\}$. Let $f_0 \colon X_0 \hra
Y=\C^2\bs \Sigma$ be the inclusion map. By the choice of $\Sigma$
there is an open tube $U\subset \C^2$ around $\Lambda$ such that
$U\cap \Sigma=\{p,q\}$. Clearly we can extend $f_0$ to a smooth
map $\tilde f\colon X\to U\bs \{p,q\} \subset Y$. However,
$f_0$ has no holomorphic extension $f\colon X\to Y$ which is
seen as follows. Any such extension would have rank at most one
by the choice of $\Sigma$. Since $f_0|_{X_0}$ already has rank one,
it would follow that $f(X)=f_0(X_0)=X_0$. But $X_0$ is a twice
punctured complex line and thus a hyperbolic space while $X$ is
covered by $\C^2$; hence $f$ (and therefore $f_0$) would be constant,
a contradiction.
\endpr

When the manifold $Y$ is Stein, the condition in
theorem 1.1 is almost necessary:

%
%  The converse
%
\proclaim 1.2 Proposition: {\rm (Gromov [Gro], 3.2.A.)}
Let $Y$ be a Stein manifold. Assume that for any Stein manifold $X$,
any closed complex submanifold $X_0 \subset X$ and any continuous map
$f_0\colon X\to Y$ which is \holo\ in a \nbd\ of $X_0$ there
exists a \holo\ map $f\colon X\to Y$ which agrees with $f_0$ to
the second order along $X_0$. Then $Y$ admits a spray.

Our proof in sect.\ 6 is slightly different from Gromov's proof.

We now turn our attention to the extension problem
for sections of \holo\ submersions. Let $h\colon Z\to X$ be
a \holo\ submersion of a complex manifold $Z$ onto a complex manifold
$X$. This means that that for each point $z\in Z$ the derivative
$d_z h$ maps $T_z Z$ surjectively onto $T_x X$, $x=h(z)$.
Let $Z_x=h^{-1}(x)$ for $x\in X$. We denote by $VT(Z)$ the
kernel of $dh$ and call it the {\it vertical tangent bundle}
(with respect to $h$). Clearly $VT(Z)$ is \holo\
subbundle of $TZ$ whose fiber at $z\in Z$ equals $$
    VT_z(Z) = \{ v\in T_z Z\colon d_zh(v) =0\} =
      T_z Z_{h(z)}.
$$

A submersion $h\colon Z\to X$ is {\it locally trivial\/} if each point
$x\in X$ has an open \nbd\ $U\subset X$ such that $h^{-1}(U)$ is
equivalent to  a product $U\times Y$ by a fiber preserving
biholomorphic map $\Phi\colon h^{-1}(U)\to U\times Y$; in such
case the submersion is a {\it \holo\ fiber bundle} over each
connected component of $X$.

A {\it section\/} of $h\colon Z\to X$ over a subset $U\subset X$
is a continuous map $f\colon U\to Z$ such that $h(f(x))=x$ for all $x\in U$.
If $U$ is an open subset of $X$ and $X_0$ is an analytic subset of
$X$, we say that a section $f$ is holomorphic on $X_0\cup U$ if
$f|_U$ is \holo\ in $U$ and $f|_{X_0}$ is holomorphic on $X_0$
(in the induced complex structure on $X_0$).

We first give a local extension theorem for sections of submersions.

\proclaim 1.3 Proposition: Let $h\colon Z\to X$ be a \holo\
submersion of a complex manifold $Z$ onto a complex manifold $X$.
Given a closed Stein subspace $X_0$ of $X$ (possibly with
singularities) and a \holo\ section $f_0\colon X_0\to Z|_{X_0}$
of $h$ over $X_0$, there are an open set $U\subset X$ containing $X_0$
and a \holo\ section $f\colon U\to Z|_U$ with $f|_{X_0}=f_0$.
If $f_0$ extends to a continuous section $f_0\colon X\to Z$,
there is a homotopy $f_t\colon X\to Z$ ($t\in [0,1]$)
which is fixed on $X_0$ such that the section
$f_1$ is \holo\ in an open set containing $X_0$.

Our proof of prop.\ 1.3 in sect.\ 6 applies also in the
case when $X_0$ is a Stein subspace of a complex space $X$,
possibly with singularities. The extension problem for maps
$f\colon X_0\to Y$ can be reduced to that of sections of the
trivial submersion $Z=X\times Y\to X$, and hence proposition 1.3
implies the existence of local \holo\ extensions in theorem 1.1.
Proposition 1.3 remains valid for families of \holo\ sections
depending continuously on a parameter in a compact Hausdorff space
$P$ (compare def.\ 3 below).

We now consider the global extension problem for sections of
submersions. If $p\colon E\to Z$ is a \hvb\ and $z\in Z$,
we write $E_z=p^{-1}(z) \subset E$ and denote by $0_z \in E_z$
the zero element of $E_z$.

%
%
%  DEFINITION OF A SPRAY
%
%
\medskip\ni \bf Definition 2.
{\rm (Gromov [Gro], sec.\ 1.1.B.)}  \sl Let  $h\colon Z\to X$ be a
\holo\ submersion and $U\subset X$ an open subset. A {\bf spray}
on $Z|_U=h^{-1}(U)$ associated to $h$ (a fiber-dominating spray in
Gromov's terminology) is a triple $(E,p,s)$, where $p\colon E\to
Z|_U$ is a \hvb\ and $s\colon E\to Z|_U$ is a \holo\ map satisfying
for each $z\in Z|_U$
\item{(i)}   $s(E_z) \subset Z_{h(z)}$
(equivalently, $h\circ p=h\circ s$),
\item{(ii)}  $s(0_z)=z$, and
\item{(iii)} the restriction of the derivative
$ds \colon T_{0_z} E \to VT_z(Z)$ to the subspace
$E_z \subset T_{0_z}E$ maps $E_z$ surjectively onto $VT_z(Z)$.
\medskip\rm

Thus a spray on a complex manifold $Y$ in the sense of definition
1 coincides with a fiber-spray associated to the trivial
submersion of $Y$ to a point. We give examples of submersions
with sprays after corollary 1.5 below.

We will have to consider parametrized families of sections
and we now introduce the relevant notions; the reader may observe
a close similarity with the objects that were considered
by Grauert [Gr2] and Cartan [Car].

%
%  P-sections
%
\medskip\ni \bf Definition 3. \sl
Let $P$ be nonempty compact Hausdorff spaces and $P_0\subset P$ a
closed subset (possibly empty) which is a strong deformation
retraction of an open \nbd\ $U\supset P_0$ in $P$.
(In our constructions $P$ will be a polyhedron and
$P_0\subset P$ a subpolyhedron.)
\item{(a)} A $P$-section of a submersion $h\colon Z\to X$ is a
continuous map $f\colon X\times  P\to Z$ such that $f_p=
f(\cdotp,p)\colon X\to Z$ is a section of $h$ for each fixed $p\in
P$. A $P$-section $f$ is holomorphic on $X$ (resp.\ on a subset
$X_0\subset X$) if $f_p$ is holomorphic on $X$ (resp.\ on $X_0$)
for each fixed $p\in P$.

\item{(b)} A homotopy of $P$-sections is a $P\times[0,1]$-section,
i.e., a continuous map $H\colon X\times P\times [0,1]\to Z$ such
that $H_t=H(\cdotp,\cdotp,t) \colon X\times P\to Z$ is a
$P$-section for each $t\in [0,1]$. The homotopy $H$ is holomorphic
if the section $H_{p,t}=H(\cdotp,p,t) \colon X\to Z$ is \holo\ for
each $(p,t)\in P\times [0,1]$.

\item{(c)} A $(P,P_0)$-section of $h$ is a $P$-section
$f\colon X\times  P\to Z$ such that the section
$f_p=f(\cdotp,p)\colon X\to Z$ is \holo\ on $X$ for each
$p\in P_0$.
\medskip\rm

The following is the main result in this paper.
%
%
%   MAIN THEOREM
%
\medskip\ni \bf  1.4 Theorem. \sl
Let $h\colon Z\to X$ be a holomorphic submersion of a complex
manifold $Z$ onto a Stein manifold $X$. Let $X_0\subset X$ be a
closed complex subvariety of $X$, $K\ss X$ a compact \hc{X}\ subset
and $U\subset X$ an open set containing $K$. Assume that for each
point $x\in X\bs K$ there is an open \nbd\ $U_x\subset X$ such
that $h\colon h^{-1}(U_x) \to U_x$ admits a spray (def.\ 2). Let
$P$ be a compact Hausdorff space. For any $P$-section
$f_0 \colon X\times P\to Z$ of $h$ which is holomorphic on $X_0\cup U$
there are an open set $U'\subset X$, with $K\subset U'\subset U$,
and a homotopy of $P$-sections $F\colon  X\times P\times [0,1]\to Z$
such that, writing $f_t=F(\cdotp,\cdotp,t) \colon X\times P\to Z$
for $t\in [0,1]$, we have
\item{(i)}   $f_0$ is the given initial $P$-section,
\item{(ii)}  the $P$-section $f_1$ is holomorphic on $X$,
\item{(iii)} For each $t\in [0,1]$ the $P$-section
$f_t$ is holomorphic on $X_0 \cup U'$, $f_t|_{X_0} = f_0|_{X_0}$,
and $f_t$ approximates $f_0$ uniformly on $K$.

\ni  Moreover, if $f_0$ is a $(P,P_0)$-section, we can we can
choose the homotopy $F$ as above such that it is fixed on $P_0$,
i.e., the section $f_{p,t}$ is independent of $t\in [0,1]$ when
$p\in P_0$.

If the initial $P$-section $f_0$ is \holo\ in an open set
$V\supset X_0\cup K$, then for any integer $k\in \Z_+$ we can
choose $F$ as above which is in addition \holo\ in a \nbd\
of $X_0\cup K$ and such that for each $(p,t)\in P\times [0,1]$,
the section $f_{p,t}$ matches $f_{p,0}$ to order $k$ along $X_0$.
(In this case it suffices to assume that the submersion $h$ admits
a spray in a small \nbd\ of each point $x\in X\bs (X_0\cup K)$,
but we do not need a spray over points in $X_0$.)
\medskip\rm

It is possible to extend theorem 1.4 to
{\it submersions with stratified sprays over
Stein spaces} (section 7).
A special case of theorem 1.4 was proved in [FP2]
(theorems 1.7 and 1.9), but the proof given there
does not carry over to the present situation.
Theorem 1.4 implies the following (see the proof of
corollary 1.5 in [FP1]):

\proclaim 1.5 Corollary: Let $h\colon Z\to X$ be a \holo\
submersion, $X_0\subset X$ a closed complex subvariety and
$f_0\colon X_0\to Z$ a \holo\ section of $h$ over $X_0$. Denote by
${\cH}(X,Z;X_0,f_0)$ (resp.\ ${\cC}(X,Z;X_0,f_0)$ the set of
global \holo\ (resp.\ continuous) sections of $h\colon Z\to X$
whose restriction to $X_0$ equals $f_0$. (Both spaces are endowed
with the compact-open topology.) Then the inclusion
${\cH}(X,Z;X_0,f_0) \hra {\cC}(X,Z;X_0,f_0)$ is a weak homotopy
equivalence, i.e., it induces an isomorphism of all homotopy
groups of the two spaces. The same is true for the inclusion
${\cH}(X,Z)\hra {\cC}(X,Z)\cap {\cH}(X_0,Z)$ of the space of
holomorphic sections into the space of continuous sections whose
restrictions to $X_0$ are holomorphic.

\demo Example 2: Fix an integer $q\ge 2$ and set
$$
    \Gamma = \{(z',z_q)\in \C^q \colon |z_q|\le 1+|z'|\}.
$$
Let $h\colon V\to X$ be a \holo\ fiber bundle with fiber
$\C^q$ and structure group $\Aut\C^q$. Let $\Sigma \subset V$
be a closed complex subvariety in $V$ satisfying:

\item{(a)} $\Sigma_x=\Sigma \cap h^{-1}(x)$ is of
complex dimension at most $q-2$ for each $x\in X$, and

\item{(b)} each point $x\in X$ has an open \nbd\ $U\subset X$
and a fiber preserving biholomorphic map
$\Phi\colon U'= h^{-1}(U) \to U\times\C^q$
such that $\Phi(U'\cap \Sigma) \subset U\times \Gamma$.

\ni Then the submersion $h\colon Z=V\bs \Sigma \to X$
admits a spray over each set $U$ as in (b)
([FP2], lemma 7.1) and hence theorem 1.4 applies.
(For \hvb s  $V\to X$ this was proved in [FP2], theorem 1.7.)
The total space $V$ of such a fiber bundle need not be
Stein even when $X=\C$ and $q=2$ [Dem].
\endpr

\demo Problem: For which compact sets $K\subset \C^n$ does the
homotopy principle hold for maps of Stein manifolds into $\C^n\bs
K$, in the sense that each continuous map $f_0 \colon X\to \C^n\bs
K$ is homotopic to a \holo\ map $f_1\colon X\to \C^n\bs K$ ?
\medskip\rm

At this point the only known examples are the finite sets.
Using analytic continuation it is easily seen that,
unless $K$ is finite, there is no spray on $Y=\C^n\bs K$
defined on a trivial bundle over $Y$ (but we don't know whether
a spray may exist on a more general \hvb s over $Y$).
For instance, does the above h-principle hold when $K$ is a
closed ball in $\C^n$ ? If so, this would give many holomorphic
maps $f\colon X\to\C^n$ with $\inf_{x\in X} |f(x)|>0$, and such
estimates are useful in the embedding-interpolation problems
as is clear from [Pre]. A good test case might be the \holo\ map
$$ {\rm SL}(2,\C)\to\C^2\bs\{0\}, \qquad
   \pmatrix{\a&\b\cr \g&\d} \to (\a,\b) \qquad  (\a\d-\b\g=1).
$$
This map is homotopic to a smooth map into the complement
of $B=\{z\in \C^n\colon |z|\le 1\}$, but it is not clear whether
it is homotopic to a \holo\ map to $\C^n\bs B$.
\endpr

We prove propositions 1.2, 1.3 and theorem 1.4 in section 6 below.
In the proof of theorem 1.4 we need the tools developed
in sect.\ 2--5. Theorem 1.1 is a special case of theorem 1.4
and will not be treated separately because we do not know
any substantially simpler proof for this special case.
In section 7 we discuss an extension of theorem 1.4
to {\it submersions with stratified sprays over Stein spaces}.

The basic analytic constructions in the proof of theorem 1.4
are similar to those in [Gra2], [Car], [Gro], [HL2], [FP1] and [FP2].
We first show how to patch holomorphic sections $a$ resp.\ $b$, defined on
sets $A\subset X$ resp.\ $B\subset X$ and extending a given \holo\
section $f_0 \colon X_0\to Z$ on a subvariety $X_0\subset X$, into
a single \holo\ section over $\tilde a$ on $A\cup B$ such that
$\tilde a|_A$ approximates $a$ and $\tilde a=f_0$ on
$(A\cup B) \cap X_0$. This can be done if $a$ and $b$ are sufficiently
close on $A\cap B$ and if the sets $A$ and $B$ satisfy certain
conditions (section 5). In the globalization scheme the
basic analytic constructions are needed for parametrized families
of sections, and for this reason we do everything by smooth
Banach space operators. We globalize the construction using the
scheme outlined in [Gro] and developed in [FP2]. We also need
certain special coverings of the base Stein manifold
constructed by Henkin and Leiterer in [HL2].

%
%
%   The h-Runge theorem
%
%
\beginsection 2. Oka-Weil theorems for submersions with sprays.

In this section we prove a homotopy version of the Oka-Weil
approximation theorem for sections of submersions with sprays
over Stein manifolds, with interpolation on a complex subvariety
$X_0$ of $X$. This is similar to results in [Gr1], [Gro]
and [FP1]. We first consider the case without parameters;
theorem 2.1 coincides with theorem 4.1 in [FP1]
when $X_0=\emptyset$.

%
%  Oka-Weil theorem, no parameters.
%
\medskip\ni \bf 2.1 Theorem. \sl
Let $h\colon Z\to X$ be a \holo\ submersion of a complex
manifold $Z$ onto a Stein manifold $X$, let $X_0 \subset X$ be a
closed complex subvariety of $X$, and let $K$ be a compact
\hc{X}\ subset. Assume that $U\supset K$ is an open set
and $f_t\colon U\to Z$ $(t\in [0,1])$ is a homotopy of holomorphic
sections of $h$ over $U$ such that $f_0$ extends to a holomorphic
section over $X$ and the homotopy is fixed over $X_0$, i.e.,
$f_t(x)=f_0(x)$ for all $x\in X_0$ and $t\in [0,1]$.
Let $d$ be a metric on $Z$ compatible with the manifold topology.
If the submersion $h$ admits a spray (def.\ 2) then for each
$\e>0$ there exists a continuous family of holomorphic
sections $\tilde f_t\colon X\to Z$ $(t\in [0,1])$ such that
\item{(a)} $\tilde f_0=f_0$,
\item{(b)} $\tilde f_t|_{X_0}=f_0|_{X_0}$ for each $t\in [0,1]$, and
\item{(c)} $d\bigl(\tilde f_t(x),f_t(x)\bigr) < \e$ for each $x\in K$ and
$t\in [0,1]$.
\medskip\rm

\demo Proof:
% This is similar to the proof of theorem 4.1 in [FP1].
Write $\wt{X}= f_0(X)\subset Z$ and
$\wt{X}_0=f_0(X_0) \subset \wt{X}$. The restriction of the spray
map $s\colon E\to Z$ to the \hvb\ $\wt{E}= E|_{\wt{X}} \to \wt{X}$
is a submersion from an open \nbd\ of the zero section
in $\wt{E}$ onto an open \nbd\ of $\wt{X}$ in $Z$.
Hence there is a $t_1>0$ such that, after shrinking $U$ slightly
around $K$, we can pull back the \holo\ sections $f_t\colon U\to Z$
for $0\le t\le t_1$ to a continuous family of \holo\
sections $\xi_t \colon f_0(U)\to \wt{E}$ such that
$\xi_0$ is the zero section, and such that $\xi_t$ vanishes
on $f_0(U\cap X_0)$ for each $t\in [0,t_1]$.

We may assume that $\bar U$ is compact. By the Oka-Cartan theory there
exist finitely many \holo\ functions $g_1,\ldots,g_k \colon X\to \C$
which vanish on $X_0$ and which generate the ideal of $X_0$ at
each point $x\in \bar U$. Since $\xi_t\circ f_0$ vanishes on
$X_0 \cap U$, we have
$$ \xi_t\bigl(f_0(x)\bigr) = \sum_{j=1}^k g_j(x)\,
   \xi^j_t\bigl(f_0(x)\bigr) \qquad (x\in U)
$$
for some holomorphic sections $\xi^j_t \colon f_0(U)\to \wt{E}$
depending continuously on $t\in [0,t_1]$. We now apply the
Oka-Weil theorem [H\"or, Theorem 5.6.2] to approximate
the family $\xi^j_t$, uniformly on $f_0(K)$, by a family of
global \holo\ sections $\tilde \xi^j_t \colon \wt{X}\to \wt{E}$,
and we set $\tilde \xi_t(z)=\sum_{j=1}^k g_j(h(z)) \tilde \xi^j_t(z)$.
The sections
$$ \tilde f_t(x)=s\bigl( \tilde \xi_t(f_0(x)) \bigr) \in Z
   \qquad (x\in X,\ t\in [0,t_1])
$$
then satisfy theorem 2.1 for $t\in [0,t_1]$.

Using $\tilde f_{t_1}$ as the new initial global section
and repeating the above construction, we obtain a $t_2>t_1$
and a family of approximating sections $\tilde f_t\colon X\to Z$
for $t\in [t_1,t_2]$. We can see as in theorem 4.1 in [FP1]
that the proof can be completed in a finite number of steps,
their number depending only on the initial family $f_t$.
\endpr

We now state the analogous approximation result for families
of sections. If $(P,P_0)$ is a pair of compact
Hausdorff spaces as in definition 3 (sect.\ 1), we set
$$
   \wt{P} = P\times [0,1],\qquad
   \wt{P}_0= (P\times \{0\})\cup (P_0\times [0,1]). \eqno(2.1)
$$
For each $\wt{P}$-section $f \colon X\times P\times [0,1]\to Z$
(def.\ 3) we write $f_t=f(\cdotp,\cdotp,t) \colon X\times P\to Z$ and
$f_{p,t}=f(\cdotp,p,t) \colon X\to Z$.

%
%
%  Oka-Weil theorem with parameters.
%
%
\medskip\ni \bf 2.2 Theorem. \sl
Let $h\colon Z\to X$, $X_0$ and $K$ be as in theorem 2.1, and
let $U,V$ be open subsets of $X$ with $K\subset U\subset V \ss X$.
Assume that $f\colon U\times \wt{P}\to Z$ is a \holo\
$\wt{P}$-section (def.\ 3) such that

\item{(i)} $f_{p,t}$ extends to a global holomorphic section on $X$
for all $(p,t)\in \wt{P}_0$ (2.1), and

\item{(ii)} $f_{p,t}(x)=f_{p,0}(x)$ for all $x\in X_0$ and
$(p,t)\in \wt{P}$. %; i.e., the homotopy is fixed on $X_0$.

\ni Let $d$ be a metric on $Z$ compatible with the manifold
topology. If the submersion $h\colon Z\to X$ admits a globally defined
spray (def.\ 2), then for each $\e>0$ there exist a \nbd\ $U'\supset K$
of $K$ and a \holo\ $\wt{P}\times [0,1]=P\times [0,1]^2$-section
$g\colon U'\times P\times [0,1]^2 \to Z$ such that,
writing $g^u=g(\cdotp,\cdotp,\cdotp,u)$ and
$g^u_{p,t}=g(\cdotp,p,t,u)\colon U'\to Z$, we have

\item{(a)} $g^0=f$, % and $g^u$ is \holo\ on $U'$ for each $u\in [0,1]$,

\item{(b)} $g^1$ extends to a \holo\ $\wt{P}$-section
$\tilde f \colon V \times \wt{P}\to Z$ over $V$,

\item{(c)} $g^u_{p,t}$ is independent of $u\in [0,1]$
when $(p,t)\in \wt{P}_0$,

\item{(d)} $g^u_{p,t}(x)=f_{p,0}(x)$ for all $x\in X_0$ and
$(p,t,u)\in P\times [0,1]^2$, and

\item{(e)} $d\bigl( g^u_{p.t}(x), f_{p,0}(x) \bigr) <\e$ for all
$x\in K$ and $(p,t,u)\in P\times [0,1]^2$.
\medskip\rm

Theorem 2.2 is proved by following the proof of theorem 4.2
in [FP1] with the obvious modifications, indicated in the
proof of theorem 2.1 above, to insure that we keep everything
fixed on the subvariety $X_0$. We omit the obvious details.

\beginsection 3. Splitting of \holo\ functions on Cartan pairs.

The main results in this section are lemmas 3.2 and 3.3 which
are needed in sect.\ 4. We shall use the following notation.
${\cH}(X)$ denotes the Fr\'echet algebra of holomorphic
functions on a complex manifold $X$, equipped with the topology
of uniform convergence on compact sets. Any closed complex
subvariety $X_0\subset X$ carries an induced structure of a reduced
complex space, and the space ${\cH}(X_0)$ of all \holo\ functions
on $X_0$ is also a Fr\'echet space (see [GuR, p.158, Theorem 5]
or [H\"or], Corollary 7.2.6). If $Y$ is another complex manifold,
${\cH}(X,Y)$ denotes the space of \holo\ maps $X\to Y$.

If $D$ is a domain in a complex manifold $X$ and $X_0$ is a closed complex
subvariety of $X$, we denote by $H^\infty_{X_0}(D)$ the Banach algebra of
bounded \holo\ functions in $D$ which vanish on $X_0\cap D$.
If $X_0\cap D=\emptyset$ we have $H^\infty_{X_0}(D)=H^\infty(D)$.
By $H^\infty(X_0\cap D)$ we denote the space of bounded \holo\
functions on the subvariety $X_0\cap D$.

%
%  Bounded extension operator.
%
\proclaim 3.1 Lemma: Let $X$ be a Stein manifold, $X_0\subset X$
a closed complex subvariety and $D\subset X$ a \psc\ domain
in $X$. Then for any relatively compact subdomain
$\Omega\ss D$ there exists a bounded linear extension
operator $S\colon H^\infty(X_0 \cap D) \to H^\infty(\Omega)$
such that $(Sf)(x)=x$ for each $f\in  H^\infty(X_0 \cap D)$
and $x\in X_0\cap \Omega$.

\demo Remark: If $D \ss X$ is \spsc\ and if $X_0$ has no
singularities on $X_0\cap bD$ and intersects $bD$ transversely,
Henkin ([Hen], [HL1]) constructed a bounded extension operator
$S\colon H^\infty(X_0 \cap D) \to H^\infty(D)$
(no shrinking of the domain!); for recent related results
see [AAC] and the references therein. No such extension
exists in general if $X_0$ has singularities along $X_0\cap bD$.

\demo Proof: We owe the idea of this proof to Bo Berndtsson
(private communication). Since $D$ is \psc\ in $X$,
the restriction operator $R \colon {\cal H}(D) \to {\cal H}(X_0\cap D)$
is surjective  [GuR, p.245, Theorem 18].
Since both spaces are Fr\'echet (and hence complete),
the open mapping theorem applies. Choose a domain $\Omega_1 \subset X$
such that $\Omega\ss \Omega_1\ss D$. By the open mapping theorem the
image by $R$ of the set
$\{f\in {\cal H}(D) \colon ||f||_{L^\infty(\Omega_1)} <1\}$
contains a \nbd\ of the origin in ${\cal H}(X_0\cap D)$.
This means that there are a relatively compact subset
$Y\ss X_0\cap D$ and a constant $M<\infty$
such that any $h \in {\cal H}(X_0\cap D)$ extends to
a function $h'\in {\cal H}(D)$ satisfying the estimate
$$ ||h'||_{L^\infty(\Omega_1)} \le M ||h||_{L^\infty(Y)}.$$
We may assume that $\Omega_1\cap X_0 \subset Y$.
The restriction $h'|_{\Omega_1}$ is bounded and hence
belongs to the Bergman space
$H= L^2(\Omega_1)\cap {\cal H}(\Omega_1)$, where the $L^2$-norm
is measured with respect to some smooth hermitian metric on $X$.
$H$ is a Hilbert space containing the closed
subspace $H_0=\{f\in H\colon f|_{X_0}=0\}$.
Let $H_1$ be the orthogonal complement to $H_0$ in $H$.
Projecting $h'$ orthogonally into $H_1$ we get a
function $\tilde h\in H_1$ which extends $h$
and which has the minimal $L^2(\Omega_1)$ norm among all
$L^2$-holomorphic extensions of $h$ to $\Omega_1$. Clearly such
$\tilde h$ is unique and $S\colon h\to \tilde h$ gives a
bounded linear operator
$S\colon H^\infty(X_0\cap D) \to L^2(\Omega_1)$.
Furthermore, restricting $\tilde h$ to the subdomain
$\Omega\ss \Omega_1$ and applying the Cauchy estimates,
we get a bounded linear extension operator
$S\colon H^\infty(X_0\cap D)\to H^\infty(\Omega)$.
\endpr

The following lemma is crucial for the results in sections 4 and 5.
If $X_0\cap \Omega=\emptyset$, lemma 3.2 coincides with
lemma 2.4 in [FP1].

%
%  A splitting lemma.
%
\medskip\ni\bf  3.2 Lemma. \sl
Let $X$ be a Stein manifold, $X_0$ be a closed complex subvariety of
$X$ and $A, B \ss X$ relatively compact domains satisfying the
following:
\item{(i)}   $\Omega=A\cup B$ is a smooth \spsc\ domain in $X$,
\item{(ii)}  $\bar{A\bs B} \cap \bar{B\bs A} =\emptyset$, and
\item{(iii)} $\bar{X_0\cap C} \subset \Omega$, where $C=A\cap B$.

\ni Then there exist bounded linear operators
${\cal A}\colon H^\infty_{X_0}(C) \to H^\infty_{X_0}(A)$,
${\cal B}\colon H^\infty_{X_0}(C) \to H^\infty_{X_0}(B)$,
such that $c={\cal A}(c) - {\cal B}(c)$ on $C$ for each
$c\in H^\infty_{X_0}(C)$.
\medskip\rm

\demo Proof: Condition (iii) implies that
$X_0 \cap b\Omega \cap \bar C=\emptyset$. Hence we
can enlarge $\Omega$ slightly around $X_0\cap b\Omega$ to
get a \spsc\ domain $\Omega'\supset \Omega$ in $X$ satisfying
$\Omega'\cap bC =\Omega\cap bC$ and
$\bar{X_0\cap \Omega} \subset X_0\cap \Omega'$. We can write
$\Omega'=A'\cup B'$ where $A'\cap B'=C$, $A'\supset A$,
$B'\supset B$, and the sets $A'$, $B'$ also satisfy
the separation property (ii).
Choose a smooth function $\chi\colon X\to [0,1]$ such that $\chi=0$
in an open \nbd\ of $\bar{A'\bs B'}$ and $\chi=1$ in an open \nbd\
of $\bar{B'\bs A'}$. For any $c\in H^\infty(C)$ the function $\chi c$
extends to a bounded smooth function on $A'$ which equals zero outside
of $C$, and likewise $(\chi-1)c$ extends to a bounded smooth function
on $B'$ which equals zero outside of $C$. The difference of these two
functions equals $c$ on $C$, but the functions are not holomorphic.

Since $\Omega'$ is a relatively compact \spsc\ domain in a Stein manifold,
there exists a linear solution operator $T$ for the $\dibar$-equation
in $\Omega'$ which is bounded in the sup-norm, i.e., for any bounded
$\dibar$-closed $(0,1)$-form $g$ on $\Omega'$ we have $\dibar(Tg)=g$ and
$||Tg||_\infty \le {\rm const}||g||_\infty$ ([HL1], p.\ 82).
Since $\supp (\dibar \chi)\cap \Omega' \subset C$, the bounded
$(0,1)$-form $g=\dibar(\chi c)=\dibar((\chi-1)c)= c\dibar \chi$
on $C$ extends to a bounded $(0,1)$-form on $\Omega'$ which is zero
outside of $C$. Set
$$ a'=\left( \chi c -T(g) \right)|_{A'}, \qquad
   b' = \left((\chi-1)c - T(g) \right)|_{B'}.
$$
It is immediate that $a'\in H^\infty(A')$, $b'\in H^\infty(B')$
and $(a'-b')|_C=c$. This solves the
problem if $X_0\cap \Omega=\emptyset$.

Suppose now that $c\in H^\infty_{X_0}(C)$. The functions $a'$ and
$b'$ need not vanish on $X_0$. However, since
$(a'-b')|_{X_0 \cap C}=c|_{X_0 \cap C}=0$, $a'$ and $b'$
define a function $h\in H^\infty(X_0\cap \Omega')$. There exists a
\psc\ domain $D\ss X$ containing $\Omega'$ such that
$\bar\Omega\subset D$ and $X_0\cap D=X_0\cap \Omega'$. Let
$S\colon H^\infty(X_0 \cap D) \to H^\infty(\Omega)$ be a bounded
linear extension operator provided by lemma 3.1.
The pair of functions
$$ a=(a'- Sh)|_A \in H^\infty_{X_0}(A),
   \qquad b=(b'-Sh)|_B \in H^\infty_{X_0}(B)
$$
then satisfies lemma 3.2, and every step in the construction was
performed by a bounded linear operator between suitable function spaces.
\endpr

\demo Remark: Lemma 3.2 also holds in spaces of
bounded \holo\ functions on $A,B,C$ which vanish to a fixed
order $k\in \Z_+$ along the subvariety $X_0$; this could
be used to give an alternative proof of the last statement
in theorem 1.4  avoiding lemma 3.3. Here is an outline
of proof. Let ${\cO}={\cO}_X$ denote the sheaf of germs
of \holo\ functions on $X$, and let ${\cJ}_k \subset {\cO}$ be
the sub-sheaf of ideals consisting of germs that vanish
to order $k$ along $X_0$. We then have
a short exact sequence of coherent analytic sheaves on $X$:
$$ 0\to {\cJ}_k \hra {\cO} \to {\cK} \to 0, $$
where the quotient sheaf $\cK={\cO}/{\cJ}_k$ is trivial on $X\bs X_0$.
Cartan's Theorem B implies that, over any Stein subset $D\subset X$,
we can lift any \holo\ section of ${\cK}$ to a \holo\ section
of ${\cO}$. If $c\in H^\infty(C)$ is a section of ${\cJ}_k$ over $C$,
we get $c=a'-b'$ with $a'\in H^\infty(A')$ and $b'\in H^\infty(B')$.
The pair $(a',b')$ defines a \holo\ section $h$ of $\cK$ over
$\Omega'=A'\cup B'$ which we can lift to a section $\tilde h$ of $\cO$
over $\Omega'$. Then the sections $a=a'-\tilde h$, $b=b'-\tilde h$
solve the problem on $(A,B)$, provided that $h \to \tilde h$ is
given by a bounded extension operator as in lemma 3.1. This can be
done exactly as before since the space of sections of the sheaf
${\cJ}_k$ (and hence of $\cK$) over $\Omega'$ is a Fr\'echet
space in the topology of uniform convergence on compact
sets ([GuR], Chapter 8).
\endpr

The next lemma is analogous to lemma 3.2, except that the
set $A\subset X$ is not relatively compact. Choose a hermitian metric
on $X$ and denote by $d\l$ the associated volume element.
For any domain $\Omega\subset X$ and
\psh\ function $\rho\colon \Omega\to \R$ let
$$ H^2_\rho(\Omega)= \{f\in {\cH}(\Omega) \colon \;
   ||f||^2_{L^2_\rho(\Omega)}=
   \int_\Omega |f|^2 e^{-\rho} d\lambda <\infty\}
$$
denote the Bergman space on $\Omega$ with weight $e^{-\rho}$.
For any subdomain $D\subset \Omega$ we define
$$
   H^{2,\infty}_{\rho}(\Omega,D)=\{f\in {\cH}(\Omega) \colon
   ||f||_{\rho,\infty} =
   ||f||_{L^2_\rho(\Omega)} + ||f||_{L^\infty(D)} < +\infty \}.
$$
Clearly this is a Banach space with the norm $||\cdotp||_{\rho,\infty}$.

%
%  Splitting lemma: unbounded domains
%
\medskip\ni \bf 3.3 Lemma: \sl
Let $X$ be a Stein manifold and $A, B \subset X$
open sets such that $\bar B$ is compact. Let $C=A\cap B$.
Assume that $\Omega=A\cup B$ is a smooth \spsc\ domain in $X$
and $\bar{A\bs B} \cap \bar{B\bs A} =\emptyset$.
For any open subset $D\ss \Omega$ there are a smooth \psh\ function
$\rho\colon \Omega \to \R$ and bounded linear operators
$$ {\cal A}\colon H^\infty(C) \to  H^{2,\infty}_{\rho}(A,A\cap D),
   \qquad
   {\cal B}\colon H^\infty(C) \to H^\infty(B)
$$
satisfying $c={\cal A}(c) - {\cal B}(c)$ on $C$ for each
$c\in H^\infty(C)$.
\medskip\rm

\demo Proof: Since we may enlarge $D\subset \Omega$ without
affecting the statement of lemma 3.3, we may assume that
$D=D'\cap \Omega$ for some domain $D'\ss X$ satisfying
$\bar B\subset D'$. Choose domains $D'_0,D'_1\ss X$
such that $D'\ss D'_0\ss D'_1$ and set
$D_j= \Omega\cap D'_j \subset\Omega$ for $j=0,1$.
We may assume that $D_1$ is \spsc. Let $T$ be a linear,
sup-norm bounded solution operator for the $\dibar$-equation
in $D_1$ ([HL1], p.\ 82). Choose a cut-off function
$\chi \colon X\to [0,1]$ as in the proof of lemma 3.2 and take
$$ \eqalign{  a'&=
   \left( \chi c -T(c\dibar\chi) \right)|_{A\cap D_1}
   \in  H^\infty(A\cap D_1), \cr
   b' &= \left((\chi-1)c - T(c\dibar \chi) \right)|_B
   \in H^\infty(B). \cr}
$$
We have $(a'-b')|_C=c$, but $a'$ is only defined on $A\cap D_1$.
To correct this we shall solve another $\dibar$-equation as follows.
Since $\Omega$ is smooth \spsc, we can enlarge it slightly
within $D'_0$ to obtain a smooth \spsc\ domain $\Omega'\subset X$
satisfying $\Omega\cup\bar D \subset \Omega'$ and
$\Omega\bs D'_0= \Omega'\bs D'_0$.
Choose a smooth cut-off function $\tau\colon X\to [0,1]$ such
that $\tau=1$ in a \nbd\ of $\bar D'_0$ and $\supp \tau\ss D'_1$.
Hence $\supp(\dibar\chi)\cap\Omega' \subset D_1\bs \bar D_0$.
The $\dibar$-closed $(0,1)$-form $g=\dibar(\tau a')=a'\dibar \tau$,
defined initially on $D_1\bs \bar D_0$, extends to a bounded form
on $\Omega'$ which is zero outside $D_1\bs \bar D_0$.
By H\"ormander [H\"or] there is a smooth \psh\ function
$\rho\colon \Omega'\to [0,+\infty)$ such that the equation
$\dibar u=g$ has a smooth solution $u$ in $\Omega'$ satisfying
$||u||_{L^2_\rho(\Omega')} \le
 {\rm const} ||g||_{\infty} \le
 {\rm const}||c||_{L^\infty(C)}$.
Moreover, if we take $u$ to be the (unique) solution
with minimal $L^2_\rho(\Omega')$ norm, the composition
$c\to g \to u$ defines a bounded linear operator
$H^\infty(C)\to L^2_\rho(\Omega')$.
By a well known estimate on the interior regularity
of the $\dibar$-operator we have for any compact
subset $K\subset \Omega'$
$$ ||u||_{L^\infty(K)} \le {\rm const}
   \bigl( ||u||_{L^2_\rho(\Omega')} +
   ||\dibar u||_{L^\infty(\Omega')} \bigr).
$$
Applying this estimate on the compact subsets
$\bar B, \bar D \subset \Omega'$ we see that
$$ a= (\tau a' - u)|_A \in H^{2,\infty}_\rho(A,A\cap D),\qquad
   b=(b'-u)|_B \in H^\infty(B).
$$
Clearly $(a-b)|_C=c$, and the maps $c\to a$, $c\to b$
for $c\in H^\infty(C)$ are bounded linear operators
into the respective Banach spaces.
\endpr

%
%
%  Attaching lemma: the model case.
%
%
\beginsection 4.  Attaching lemma: the model case.

In this section we apply lemmas 3.2 and 3.3 to construct \holo\
sections of certain model fibrations. An iteration scheme to solve
this problem was proposed in [Gro]; we shall apply the
implicit function theorem in suitable function spaces.

Recall that a pair of open subset $U\subset V$ in a complex
manifold $X$ is a {\it Runge pair}, or $U$ is Runge in $V$,
if every \holo\ function in $U$ can be approximated, uniformly
on compacts in $U$, by functions holomorphic in $V$.
We use the notation introduced in sect.\ 3 above, and
we identify the subvariety $X_0\subset X$
with $X_0\times \{0\} \subset X\times \C^n$.

%
%  Gluing-model case
%
\medskip\ni\bf 4.1 Proposition.  \sl
Let $X$ be a Stein manifold, $X_0$ a closed complex subvariety of
$X$ and $A,B\ss X$ relatively compact open subsets as in
lemma 3.2. Let $\wt{B}\supset \bar B$ and
$\wt{C}\supset \bar C$ be open sets such that
$\wt{C}\subset \wt{B}$ and $\wt{C}$ is Runge
in $\wt{B}$. Let $U\subset \cn$ be an open \nbd\ of the origin
and $\psi_0\colon \wt{C}\times U\to\cn$ a bounded
\holo\ map such that for each $x\in \wt{C}$ we have
$\psi_0(x,0)=0$ and $\psi_0(x,\cdotp) \colon U\to\cn$ is injective.
Then there exist an open \nbd\ ${\cal W}$ of $\psi_0$ in the Banach
space $H^\infty_{X_0} (\wt{C}\times U)^n$ and smooth Banach
space operators
${\cal A'}\colon {\cal W}\to H^\infty_{X_0}(A)^n$,
${\cal B'}\colon {\cal W}\to H^\infty_{X_0}(B)^n$,
with ${\cal A'}(\psi_0)=0$ and ${\cal B'}(\psi_0)=0$, such that for
each $\psi\in {\cal W}$ the  bounded \holo\ maps
$\alpha={\cal A'}(\psi) \colon A\to\cn$ and
$\beta={\cal B'}(\psi) \colon  B\to\cn$ satisfy
$\alpha|_{A\cap X_0}=0$, $\beta|_{B\cap X_0}=0$ and
$$
   \psi\bigl(x, \alpha(x) \bigr)= \beta(x) \quad (x\in C=A\cap B).
                                                         \eqno(4.1)
$$
Moreover, if $\psi \in {\cal W}$ satisfies $\psi(x,0)=0$ for
$x\in\wt{C}$, then ${\cal A'}(\psi)=0$ and ${\cal B'}(\psi)=0$.
\medskip\rm

\demo Remark: We can view a pair of maps satisfying (4.1) as a
section of the fibration over $A\cup B \subset X$ obtained
by identifying the point $(x,z)\in C\times U\subset A\times \cn$
with the point $(x,\psi(x,z)) \in C\times \C^n$. For applications
to parametrized families it is convenient to have solutions
given by operators, although this could be avoided by a suitable
analogue of Satz 8 in [Gr1].

\demo Proof: We apply the proof of proposition 5.2 in [FP1],
replacing lemma 2.4 in [FP1] by lemma 3.2 above.
We recall the main idea and refer to [FP1] for the details.
Assume first $\psi_0(x,z)=z$ $(x\in \wt{C},\ z\in U)$.
Let ${\cal A}\colon H^\infty_{X_0}(C)^n \to H^\infty_{X_0}(A)^n$
and ${\cal B}\colon H^\infty_{X_0}(C)^n \to H^\infty_{X_0}(B)^n$ be
the linear operators obtained by applying lemma 3.2
componentwise; hence $c= {\cal A}c - {\cal B}c$ for
all $c\in H^\infty_{X_0}(C)^n$. Consider the operator
$$ \eqalign{ & \Phi\colon H^\infty_{X_0}(\wt{C}\times U)^n
               \times H^\infty_{X_0}(C)^n
   \to H^\infty_{X_0} (C)^n, \cr
    & \Phi(\psi,c)(x) = \psi\bigl( x,{\cal A}c(x) \bigr) - {\cal B}c(x)
    \qquad (x\in C).  \cr}
$$
This is a smooth Banach space operator in a \nbd\ of the point
$(\psi_0,0)$, satisfying $\Phi(\psi_0,c)={\cal A}c-{\cal B}c=c$.
By the implicit function theorem the equation $\Phi(\psi,c)=0$ has
locally near $(\psi_0,0)$ a unique solution $c={\cC}(\psi)$
given by a smooth operator ${\cC}$. The operators
${\cal A'} = {\cal A}\circ{\cal C}$,
${\cal B'}={\cal B}\circ {\cal C}$ then satisfy proposition 4.1.
The general case is reduced to this one by approximating
$\psi_0$ by a \holo\ map on $\wt{B}\times U$ (see [FP1]).
\endpr

Applying lemma 3.3 instead of lemma 3.2 in the proof
of proposition 4.1 we get the following result.

%
%
%
%\vfill\eject
\medskip\ni\bf 4.2 Proposition.  \sl
Let $X$ be a Stein manifold and let $A,B\subset X$ and
$D\subset \Omega=A\cup B$ be open subsets as in lemma 3.3.
Let $\wt{B}\supset \bar B$ and $\wt{C}\supset \bar C$ be
open sets such that $\wt{C}\subset \wt{B}$ and $\wt{C}$ is Runge
in $\wt{B}$. Let $U\subset \cn$ be an open \nbd\ of the origin
and $\psi_0\colon \wt{C}\times U\to\cn$ a bounded
\holo\ map such that for each $x\in \wt{C}$,
$\psi_0(x,0)=0$ and $\psi_0(x,\cdotp) \colon U\to\cn$
is injective. Then there exist an open \nbd\
${\cal W} \subset H^\infty (\wt{C}\times U)^n$
of $\psi_0$ and smooth Banach space operators
${\cal A'}\colon {\cal W}\to H^{2,\infty}_{\rho}(A,A\cap D)^n$,
${\cal B'}\colon {\cal W}\to H^\infty(B)^n$,
with ${\cal A'}(\psi_0)=0$ and ${\cal B'}(\psi_0)=0$, such
that for each $\psi\in {\cal W}$ the \holo\ maps
$\alpha={\cal A'}(\psi) \colon A\to\cn$ and
$\beta={\cal B'}(\psi) \colon  B\to\cn$ satisfy (4.1).
Moreover, if $\psi \in {\cal W}$ satisfies $\psi(x,0)=0$ for
all $x\in\wt{C}$, then ${\cal A'}(\psi)=0$ and ${\cal B'}(\psi)=0$.
\medskip\rm

%
%
%  Sect. 5: Gluing sections over Cartan pairs.
%
%
\beginsection 5. Attaching lemma for holomorphic sections of submersions.

In this sections we use the results of sect.\ 4 in order
patch together \holo\ sections of a submersion $Z\to X$
over a Cartan pair $(A,B)$ in $X$, provided that the two
sections are sufficiently close on the intersection $A\cap B$.
Theorem 5.1 below is similar to results in sect.\ 1.6
of [Gro] and to theorems 5.1 and 5.5 in [FP1]; the additional
point here is the interpolation on a subvariety $X_0 \subset X$.
Theorem 5.2 is new and depends on proposition 4.2.

We recall from [FP2] the definition of a Cartan pair.
(Note that our earlier definition in [FP1] did not include
the Runge property (iii).)

%
%  Cartan pairs
%
\medskip\ni \bf Definition 4. \sl
An ordered pair of compact sets $(A,B)$ in a complex manifold
$X$ is said to be a {\bf Cartan pair} if
\item{(i)}    each of the set $A$, $B$, and $A\cup B$ has a
basis of Stein \nbd s,
\item{(ii)}   $\bar{A\bs B} \cap \bar{B\bs A} =\emptyset$
\ \ (separation condition), and
\item{(iii)}  the set $C=A\cap B$ is Runge in $B$.
($C$ may be empty.)

\ni If in addition $X_0$ is a closed complex subvariety of $X$
such that $X_0\cap C \subset {\rm Int}(A\cup B)$, we shall
say that the Cartan pair $(A,B)$ is $X_0$-regular.
\medskip\rm

We proved in [FP1] that for each Cartan pair $(A,B)$ in $X$
there exist bases of decreasing open \nbd s $A_j\supset A$,
$B_j\supset B$ $(j\in \Z_+)$ such that each pair $(A_j,B_j)$
satisfies the hypothesis (i) and (ii) of lemma 3.2 and
$C_j$ is Runge in $B_j$. The same proof shows that, if $(A,B)$
is $X_0$-regular, we also get property (iii) in lemma 3.2
for each $(A_j,B_j)$. Hence proposition 4.1 can be applied on
a suitable bases of \nbd s of any $X_0$-regular Cartan
pair $(A,B)$ in $X$. This allows us to glue sections
of submersions $h\colon Z\to X$ over $(A,B)$ when $Z$ admits
a spray over a \nbd\ of $B$.

The presence of parameters and the need to do everything
by homotopies complicates the statement, so let us first
explain the result in the basic case.
Let $(A,B)$ be an $X_0$-regular Cartan pair in $X$.
We are given \holo\ sections $a\colon \wt{A}\to Z$
resp.\ $b\colon \wt{B}\to Z$ over open sets
$\wt{A}\supset A$ resp.\ $\wt{B}\supset B$ such that
$a$ and $b$ agree on $X_0\cap \wt{C}$, where
$\wt{C}=\wt{A}\cap \wt{B}$.
If $Z$ admits a spray over $\wt{B}$ and if $\tilde b$
is sufficiently close to $\tilde a$ on $\wt{C}$,
we can move each of the two sections a little
(by \holo\ homotopies of sections on $\wt{A}$
resp.\ $\wt{B}$) such that the final pair
of sections coincides on $\wt{A}\cap \wt{B}$
and hence gives a \holo\ section over $\wt{A}\cup\wt{B}$.
(We must shrink the \nbd s of $A$ and $B$ in the process.)
Moreover, we can perform the procedure so that the homotopy
is fixed on $X_0$ where the initial sections already agree.
The same can be done for $P$-sections so that the
homotopies are fixed for those values of the parameter
for which the two initial sections already coincide
over $\wt{C}$.

Having said this, we state the result in precise terms.

%
%
%  Gluing of sections: the general case with parameters.
%
%
\medskip\ni\bf 5.1 Theorem. \sl
Let $h\colon Z\to X$ be a holomorphic submersion onto a
Stein manifold $X$, let $X_0$ be a closed complex
subvariety of $X$, and let $(A,B)$ be a
$X_0$-regular Cartan pair in $X$ (def.\ 4). Suppose that $\wt{B}$ is
an open \nbd\ of $B$ in $X$ such that the restriction
$Z|_{\wt{B}} = h^{-1}(\wt{B})$ admits a spray over $\wt{B}$ (def.\ 2).
Let $(P,P_0)$ be as in def.\ 3. Let $\wt{A}\supset A$ be an
open \nbd\ of $A$ in $X$ and $a\colon \wt{A} \times P\to Z$ a
\holo\ $P$-section over $\wt{A}$ (def.\ 3).
Fix a metric $d$ on $Z$ compatible with the manifold topology.
Then for each $\e>0$ there is a $\d>0$ satisfying the following
property. If $b\colon \wt{B} \times P \to Z$ is a \holo\ $P$-section
over $\wt{B}$ satisfying
$$ \eqalign{ & d\bigl( a_p(x), b_p(x)\bigr) <\d \quad
    (x\in \wt{C}=\wt{A}\cap \wt{B},\ \ p\in P), \cr
    & a_p(x)=b_p(x) \quad\quad\quad\ (x\in \wt{C},\ p\in P_0)
    \ {\rm or}\ (x\in X_0\cap \wt{C},\ p\in P), \cr}
    % & a_p(x)=b_p(x) \quad\quad\quad\ (x\in X_0\cap \wt{C},\ \ p\in P),
    % \cr}
$$
then there exist smaller \nbd s $A'\supset A$, $B'\supset B$ and
homotopies of \holo\ $P$-sections $a^t \colon A'\times P\to Z$,
$b^t \colon B'\times P\to Z$ $(t\in [0,1])$
such that

\item{(a)} $a^t_p=a_p$ and $b^t_p=b_p$ for
$(p,t)\in \wt{P}_0= (P\times \{0\})\cup (P_0\times [0,1])$,

\item{(b)} $a^t_p(x)=a_p(x)$ and $b^t_p(x)=b_p(x)$ for
$x\in X_0$ and $(p,t)\in P\times [0,1]$,

\item{(c)} $a^1_p(x)=b^1_p(x)$ for $x\in C'=A'\cap B'$ and $p\in P$, and

\item{(d)} for each $(p,t)\in P\times [0,1]$ we have
$$ d\bigl( a^t_p(x), a_p(x) \bigr) < \e \ \ (x\in A'), \quad
   d\bigl( b^t_p(x), b_p(x) \bigr) < \e \ \ (x\in C').
$$
\medskip\rm

\demo Remark: Property (c) implies that $a^1$ and $b^1$ together
define a \holo\ $P$-section
$\tilde a\colon (A'\cup B')\times P\to Z$ over $A'\cup B'$.
For $p\in P_0$ the section $\tilde a_p$ agrees with $a_p$ and
$b_p$ according to (a), and the homotopy is fixed on $X_0$
according to (b).

\demo Proof: It suffices to apply the proof of theorem 5.5 in
[FP1], but replacing proposition 5.2 in [FP1] by proposition 3.2
above. We recall the main idea since we shall need this in
the next theorem. Suppose for simplicity that $P$
is a singleton and $P_0=\emptyset$.
We linearize the problem by first constructing for some
large $n\in \Z_+$ a pair of \holo\ maps
$$ s_1 \colon V\subset \wt{A} \times \C^n \to Z, \qquad
   s_2 \colon \wt{B}\times \C^n \to Z,              \eqno(5.1)
$$ where $V$ is an open set containing $\wt{A}\times \{0\}$,
such that $s_1$ and $s_2$ are fiber preserving (i.e., the fiber
over $x\in \wt{A}$ resp.\ $x\in \wt{B}$ is mapped into
$Z_x=h^{-1}(x)$), and they are submersions of open \nbd s of the
zero sections in $\wt{A}\times \C^n$ resp.\ $\wt{B}\times
\C^n$ onto open \nbd s of the graphs $a(\wt{A})\subset Z$ resp.\
$b(\wt{B})\subset Z$. Moreover we have $s_1(x,0)=a(x)$ and
$s_2(x,0)=b(x)$.

The submersion $s_1$ is only locally defined and
can be obtained from local flows of vertical \hvf s on $Z$
(tangent to fibers of $h$) near $a(A)$. The second map $s_2$ is
globally defined and is obtained by restricting the spray map
$s\colon E\to Z|_{\wt{B}}$ to the vector bundle $E|_{b(\wt{B})}$
over the section $b(\wt{B}) \subset Z$. The total space
$E|_{b(\wt{B})}$ is not necessarily trivial; however,
we may choose $\wt{B}$ to be Stein, and hence there is a
\hvb\ epimorphism $\theta \colon \wt{B}\times \C^n \to E|_{b(\wt{B})}$
for any sufficiently large $n$. The composition
$s_2=s\circ\theta$ satisfies our requirements.

When $a$ and $b$ are sufficiently close to each other over $\wt{C}$,
we can construct a holomorphic transition map
$\psi \colon \wt{C}\times U\to \C^n$ as in
proposition 4.1 such that
$$ s_1(x,z)= s_2(x,\psi(x,z)) \qquad
   (x\in \wt{C},\ z\in U \subset \cn).              \eqno(5.2)
$$
The closeness of $a$ and $b$ over $\wt{C}$ implies that
$\psi$ is close to a map $\psi_0$ which preserves zero section:
$\psi_0(x,0)=0$ for $x\in \wt{C}$.
If $\alpha\colon A'\to U\subset \C^n$ and
$\beta\colon B'\to \C^n$ are \holo\ maps on open sets
$A'\subset \wt{A}$ resp.\ $B'\subset \wt{B}$
as in proposition 4.1, satisfying
$\psi(x,\alpha(x))=\beta(x)$ for $x\in C'=A'\cap B'$,
we set for each $t\in [0,1]$
$$ a^t(x)=s_1(x,t\alpha(x))\ \ (x\in A'),\qquad
   b^t(x)=s_2(x,t\beta (x))\ \ (x\in B').              \eqno(5.3)
$$
Then $a^0=a$, $b^0=b$, and $a^1(x)= b^1(x)$ for $x\in C'$;
hence $a^1$ and $b^1$ together define a \holo\ section
$\tilde a \colon A'\cup B'\to Z$. The homotopies $a^t$ and $b^t$
are fixed on $X_0$. The details can be found in [FP1].
\endpr

We need a similar result obtained from proposition 4.2.
For simplicity we state the result without parameters
(when $P$ is a singleton), even though the result holds
in the same generality as theorem 5.1 above.

%
%  Gluing sections on unbounded sets.
%
\medskip\ni\bf Theorem 5.2. \sl
Let $h\colon Z\to X$ be a holomorphic submersion onto a
Stein manifold $X$ and let $X_0$ be a closed complex
subvariety of $X$. Furthermore let $A\subset X$, $B\subset X$
be closed subsets such that $B$ is compact while $A$
contains a \nbd\ of $X_0$. Assume that
\item{(i)}   $A\cup B$ is (the closure of) a \spsc\ domain in $X$,
\item{(ii)}  $\bar{A\bs B} \cap \bar{B\bs A} =\emptyset$, and
\item{(iii)} the set $C=A\cap B$ is Runge in $B$ and $C\cap X_0=\emptyset$.

Assume that $Z$ admits a spray over an open set $\wt{B}\supset B$.
Let $\wt{A}\supset A$ be an open set and $a\colon \wt{A} \to Z$ a
\holo\ section. Fix a metric $d$ on $Z$ compatible with the manifold
topology. Let $D\ss X$ be an arbitrary relatively compact subset
containing $\bar B$. Choose an integer $k\in \Z_+$.
Then for each $\e>0$ there is a $\d>0$
satisfying the following property.
If $b\colon \wt{B} \to Z$ is a \holo\ section satisfying
$d\bigl( a(x), b(x)\bigr) < \d$ for
$x\in \wt{C}=\wt{A}\cap \wt{B}$,
there exist open sets $A',B'\subset X$ satisfying
$X_0\cup (A\cap D) \subset A'\subset \wt{A}$,
$B\subset B'\subset \wt{B}$, and homotopies of \holo\
sections $a^t \colon A'\to Z$, $b^t \colon B'\to Z$ $(t\in [0,1])$
such that

\item{(a)} $a^0=a$ and $b^0=b$,

\item{(b)} $a^1(x)=b^1(x)$ for $x\in C'=A'\cap B'$,

\item{(c)} for each $t\in [0,1]$, $a^t$ agrees
with $a=a^0$ to order $k$ along $X_0$ and satisfies
$$ d\bigl( a^t(x), a(x) \bigr) < \e \ \ (x\in A'\cap D),
   \quad
   d\bigl( b^t(x), b(x) \bigr) < \e \ \ (x\in C').
$$
\medskip\rm

\demo Proof: The proof is essentially the same as in theorem 5.1
except for the construction of the map $s_1$ (5.1).
As in theorem 5.1 we construct a preliminary fiber-preserving
\holo\ submersion $\tilde s_1\colon V\to Z$ from an open set
$V\subset \wt{A}\times \C^n$, with $\wt{A}\times \{0\} \subset V$,
onto a \nbd\ of $a(\wt{A})$ in $Z$. By the Oka-Cartan theory there exist
finitely many \holo\ functions $h_j\colon X\to \C$ $(1\le j\le m)$ such
that $X_0=\{x\in X\colon h_j(x)=0,\ 1\le j\le m\}$ and each $h_j$
vanishes to order $k$ on $X_0$.
Let $g\colon \wt{A}\times \C^{nm}\to \C^n$ be the map
$g(x,v_1,\ldots,v_m)= \sum_{j=1}^m h_j(x) v_j$,
where $x\in\wt{A}$ and $v_j\in \C^n$ for each $j$.
Clearly $g(x,\cdotp)$ is a linear epimorphism for each $x$
outside $X_0\cap \wt{A}$ and is degenerate for $x\in X_0\cap \wt{A}$.
Set $s_1=\tilde s_1\circ (Id,g)\colon \wt{A}\times\C^{N} \to Z$,
where $Id$ indicates the identity on $\wt{A}$ and $N=nm$.
Then $s_1$ is a submersion on a \nbd\ of
$(\wt{A} \bs X_0)\times \{0\}$ in $\wt{A}\times \C^N$
(and is degenerate over $X_0\cap \wt{A}$). We construct
a map $s_2\colon \wt{B}\times \C^N\to Z$ (5.1) as before.

Since $C=A\cap B$ does not intersect $X_0$, we may
assume that the closure of $\wt{C}=\wt{A}\cap \wt{B}$
does not intersect $X_0$ either. Hence $s_1$ is a submersion
over $\wt{C}$ which allows us to construct the transition
map $\psi$ (5.2) as before. For any solution
$(\alpha,\beta)$ of the equation (4.1) we
get the corresponding pair of sections
$$ a^1(x)= s_1(x,\alpha(x))= \tilde s_1(x,g(x,\alpha(x))),\qquad
   b^1(x)= s_2(x)
$$
as in (5.3) which agree on a \nbd\ of $C$. The section
$b^1$ is defined on a \nbd\ $B'\supset B$ as before.
On the other hand, the domain of $a^1(x)$ may shrink because
we cannot control the sup-norm of $\alpha$ on all of $A$
(since $A$ is unbounded), and hence the section
$\tilde \alpha(x)= (x,g(x,\alpha(x))) \in\wt{A}\times \C^n$
may escape from the domain of $\tilde s_1$. However,
things are not too bad. We can control the sup-norm
of $\alpha$, and hence of $\tilde \alpha$, on $A\cap D$ for
any $D\ss X$. By construction the map
$x\to g(x,\alpha(x)) =\sum_j h_j(x)\alpha_j(x)$
vanishes to order $k$ on $X_0$; hence $X_0$ is in the
domain of $a^1$ as well. It follows that
$a^1$ is \holo\ on a set $A'\supset X_0\cup (A\cap D)$
and it agrees with $a=a^0$ to order  $k$ along $X_0$.
The same applies to each section in the homotopy $a^t$
from $a=a^0$ to $a^1$, and this insures the validity
of property (c). All the rest is the same as in theorem 5.1;
in particular we get uniform approximation in (d)
from the uniform estimates on $A\cap D$ in proposition 4.2.
\endpr

\beginsection 6. Proof of propositions 1.2, 1.3 and theorem 1.4.

\demo Proof of proposition 1.2:
We embed $Y$ as a closed complex submanifold
of some Euclidean space $\C^N$. By Docquier and Grauert
([GuR], p.257, Theorem 8) there is a holomorphic retraction
$\pi\colon U\to Y$ of an open \nbd\ $U\subset \C^N$ onto $Y$.
We identify the \holo\ tangent bundle $X=TY$ with a subbundle
of $T\C^N|_Y$ and we identify $Y$ with the zero section of $TY$.
Denote the points in $TY$ by $(y,\xi)$ and let $0_y=(y,0)$.
There is an open \nbd\ $V\subset TY$ of the zero section on
which the map $s_0(y,\xi) = \pi(x+\xi) \in Y$ is defined and
holomorphic. By modifying $s_0$ outside a \nbd\ of the zero section
we can extend it to a smooth map $s_0\colon TY\to Y$. Clearly the
derivative $ds_0 \colon T_{0_y}(TY) \to T_y Y$ restricts to the
identity map on $T_y Y \subset T_{0_y}(TY)$ for each $y\in Y$.
Hence $s_0$ satisfies the requirements for a spray,
except that it is not globally holomorphic on $TY$.
Observe that $X=TY$ is a Stein manifold and
$X_0=Y$ (the zero section) is a closed complex submanifold of $X$.
If $s\colon TY\to Y$ is a \holo\ map which matches $s_0$ to the
second order along the zero section (such $s$ exists by the
hypothesis in proposition 1.2), then $s$ is a spray on $Y$.
\endpr

\demo Proof of proposition 1.3:
The image $Z_0=f_0(X_0)$ is a closed Stein subspace
of $Z$ and hence it has an open Stein \nbd\ in $Z$ according to
[Siu]. The same applies to $X_0$ in $X$; hence we may assume that
both $X$ and $Z$ are Stein. For each $g \in {\cal O}(Z)$ we denote
by $d'g$ the restriction of the differential $dg$ to the
vertical bundle $VT(Z)$. By the Oka-Cartan's theory there are
functions $g_1,\ldots,g_d \in {\cO}(Z)$ such that
$g_j=0$ on $Z_0$ for each $j$ and $\{d'g_j \colon 1\le j\le d\}$
span the vertical contangent space $VT^*_z(Z)$ at each point
$z\in Z_0$. Hence the map
$$ G\colon Z\to X\times \C^d,\qquad
   G(z)= \bigl( h(z), g_1(z),\ldots,g_d(z)\bigr)
$$
embeds a \nbd\ of $Z_0$ as a closed complex submanifold $W$
in an open set $V \subset X\times \C^d$,
with $G(Z_0) = X_0\times \{0\}^d$. Denoting by
$p\colon X\times \C^d\to X$ the projection
onto $X$, we have $p\circ G=h$.

We claim that, after shrinking $V$ around $W$, we can find
a holomorphic retraction $\pi\colon V\to W$ satisfying $p\circ \pi=p$.
Such a retraction is constructed as in the Docquier--Grauert
theorem ([GuR], p.257, Theorem 8); here is a brief outline.
Let $VT(W)\to W$ be the vertical tangent bundle to $W$
(with respect to $p\colon W\to X$); observe that $VT(W)$
is a subbundle of the trivial bundle $W\times \C^d$.
Choose a complementary \hvb\ $H\to W$ such that
$H\oplus VT(W)=W\times \C^d$. Denote the points
of $H$ by $(x,w,\xi)$, with
$(x,w)\in W$ and $\xi\in H_{(x,w)}\subset \C^d$. Consider the map
$H\to X\times \C^d$, $(x,w,\xi)\to (x,w+\xi)$
(the addition takes place in the fiber $\{x\}\times \C^d$).
As in the proof of the Docquier--Grauert theorem we see that this
map takes a \nbd\ of the zero section in $H$ biholomorphically
onto a \nbd\ of $W$ in $X\times \C^d$; hence
it conjugates the base projection $H\to W$
to a desired retraction $V \to W$.

To extend the initial \holo\ section we take
$f(x)=G^{-1}(\pi(x,0))$, where $0$ denotes the
origin in $\C^d$. When $x$ is sufficiently near $X_0$,
the point $(x,0) \in X\times \C^d$ belongs to $V$
(the domain of $\pi$), and hence $\pi(x,0)\in W$ belongs
to the range of $G$. Thus $f$ is well-defined and
\holo\ in an open set $U\supset X_0$ in $X$, and we have
$f(x)\in Z_x$ for $x\in U$.
When $x\in X_0$, we have $(x,0)\in W$, hence
$\pi(x,0)=(x,0)=G(f_0(x))$ and therefore $f(x)=f_0(x)$.

If $f_0$ extends continuously to a section $X\to Z$,
we can patch the section $f$ obtained above
with $f_0$ in a small \nbd\ of $X_0$ as follows.
As above we have an embedding of a \nbd\ of $Z_0$
in $Z$ as a submanifold $W \subset X\times \C^d$;
we patch the two sections in the ambient space
$X\times \C^d$ (where we have linear fibers) by
a cut-off function, and finally we project the result
back to $W$ by the holomorphic retraction.
The new section $f_1 \colon X\to Z$ equals $f$ in
a \nbd\ of $X_0$ (so it is \holo\ there), and it
equals $f_0$ outside a larger \nbd\ of $X_0$.
The construction also gives a homotopy
of $f_0$ to $f_1$ which is fixed on $X_0$.
This proves proposition 1.3. The proof carries
over verbatim to the case when $X$ is a complex space.
\endpr

\demo Proof of theorem 1.4:
We consider two cases: in the first case the initial
continuous section $f_0 \colon X\to Z$ is holomorphic on the
subvariety $X_0$ and in a \nbd\ of a \hc{X}\ set $K\subset X$;
in the second case $f_0$ is assumed to be \holo\ in an open
\nbd\ of $X_0\cup K$. When $K=\emptyset$, we can reduce
the first case to the second one by proposition 1.3;
however, this reduction does not really simplify the proof,
its only apparent advantage being that we need not assume
the existence of a spray on $Z$ over points in $X_0$
which turns out to be very convenient in applications.

In the first case we can obtain the required homotopy
$f_t \colon X\to Z$ satisfying theorem 1.4 by following
the proof of theorem 1.5 in [FP2] (sect.\ 6), except that
we replace the approximation and patching results used there
by the corresponding results with interpolation on $X_0$,
given by theorems 2.2 and 5.1 above. The heart of the proof
is an induction scheme (sect.\ 6 in [FP2]) in which the
main ingredient is proposition 5.1 from [FP2].

However, the second case (when $f_0$ is \holo\ in a \nbd\ of
$X_0\cup K$) requires some modifications because we
patch pairs of \holo\ sections on so-called Cartan
pairs $(A,B)$ in $X$ where the set $A$ is unbounded
(it contains $X_0\cup K$). We do this by replacing
theorem 5.1 by theorem 5.2 above. In this case the
globalization process requires small modifications,
and we feel that it would be dishonest to leave this entirely
to the reader. So we shall stir a middle course by indicating
the essential steps of the argument in both cases and
referring to [FP2] for the details.

First we must recall from sect.\ 3 in [FP2] the notion of \holo\
(and continuous) complexes and prisms associated to a given open
covering $\cU=\{U_j\}$ of the base manifold $X$. Given
such a covering $\cU$, we denote by $\cK(\cU)$ its nerve
(an infinite combinatorial simplicial complex) and by $K(\cU)$
its geometric realization (an infinite polytope).

Let $h\colon Z\to X$ be a given submersion.
A {\it holomorphic $\cK(\cU)$-complex with values in $Z$} is a family
$f_*=\{f_t\colon t\in K(\cU) \}$ of holomorphic sections of $Z$,
depending continuously on the parameter $t\in K(\cU)$, where the domain of
the section $f_t \in f_*$ is determined as follows: If $t$ belongs
to the $k$-dimensional simplex in $K(\cU)$ which is determined by
$k+1$ open sets $U_{j_0},U_{j_1},\ldots,U_{j_k} \in \cU$, then
$f_t$ is a \holo\ section of $Z$ over the set $\cap_{i=0}^k U_{j_i}$.
We also have natural restriction conditions for sections
in the family $f_*$ on boundaries of simplices in $K(\cU)$ (see
[FP2]). Thus the vertices of $K(\cU)$ correspond to \holo\
sections on the sets $U_j\in \cU$, the edges correspond to
one-parameter homotopies of sections defined on the intersections
$U_i\cap U_j$, etc. A global section $f \colon X\to Z$ may be considered
as a `constant complex', meaning that any section in the associated
complex $f_*$ is the restriction of $f$ to the appropriate
open set in $X$. A {\it holomorphic $k$-prism} is a homotopy of
holomorphic complexes $f_{*,s}$ depending on a parameter $s\in [0,1]^k$.
Similarly one defines continuous complexes and prisms with values in $Z$
as collections of continuous sections. We also work with
coverings $\cA=\{A_i\}$ of $X$ consisting of compact sets;
a $\cK(\cA)$-complex is represented by $\cK(\cU)$-complexes
for open coverings $\cU=\{U_i\}$ of $X$ with $U_i\supset A_i$
for each $i$, and we identify two complexes whose sections
agree near the corresponding sets in $\cA$. For details we
refer to sect.\ 3 in [FP2].

Assume for simplicity that $P$ is a singleton; the proof
in the general case follows the same pattern.
Consider first the case when the restriction $f_0|_{X_0}$
is \holo\ on $X_0$ and $f_0$ is \holo\ over a \nbd\
$U\supset K$ and of a compact \hc{X}\ subset $K\subset X$.
To harmonize the notation with [FP2] we write $f_0=a$.
Our goal is to move $a$ by a homotopy which is fixed
on $X_0$ to a \holo\ section $f\colon X\to Z$.

By theorem 4.6 in [FP2] (which follows from the results
in sect.\ 2 of [HL2]) there exists a sequence
${\cA}=(A_0,A_1,A_2,\ldots)$ of compact,
\hc{X} subsets in $X$, satisfying

\item{(i)}   $K\subset A_0\subset U$ and $X=\cup_{j=0}^\infty A_j$,

\item{(ii)}  for each $n\in \Z_+$ the pair of sets
$(A^n,A_{n+1})$, where $A^n=A_0\cup A_1\cup\cdots\cup A_n$,
is a $X_0$-regular Cartan pair in the sense of def.\ 4 above, and

\item{(iii)} the submersion $h\colon Z\to X$ admits a
spray over an open \nbd\ of $A_j$ for each $j\ge 1$
(but we don't need a spray over $A_0$).

\ni  As in [FP2] we call any finite sequence satisfying
(ii) and (iii) a {\it Cartan string}, and the entire
collection $\cA$ is a {\it Cartan covering\/} of $X$.
The order of sets in a Cartan string is important
because $A^n \cap A_{n+1}$ must be Runge in $A_{n+1}$.
In [FP2] there was no subvariety $X_0$,
but the construction in [HL2] shows that we
can easily satisfy the $X_0$-regularity condition in (ii) by
a suitable choice of the set $A_n$ at each step.

The sets $A_j$ for $j\ge 1$ may be chosen arbitrarily small,
subordinate to any open covering of $X$. By proposition
4.7 in [FP2] we can deform the initial section $a$ (which we
now consider as a constant $\cK(\cA)$-complex $a_{*,0}$)
by a homotopy of continuous $\cK(\cA)$-complexes $a_{*,s}$
$(s\in [0,1])$ into a \holo\ $\cK(\cA)$-complex $a_{*,1}$.
The only additional requirement here is that all sections
belonging to any of the complexes $a_{*,s}$ $(s\in [0,1])$
must agree with the initial section $a$ on the intersection
of their domains with the subvariety $X_0$. This is easily
satisfied by first choosing local \holo\ extensions of $a|_{X_0}$;
the homotopies between them are obtained by taking their
convex linear combinations with respect to a local
linear structure on the fibers of $Z$ (see [FP2]),
and hence all sections will agree with $a$ on $X_0$.

Applying proposition 5.1 in [FP2] to the \holo\ complex
$a_{*,1}$ we inductively construct a sequence of \holo\
$\cK(\cA)$-complexes $a_{*,n}$ ($n=1,2,3,\ldots$) which
converges as $n\to\infty$ to a holomorphic section
$a_{*,\infty}=f\colon X\to Z$. The basic tools that we use
in proposition 5.1 from [FP2] in the present case are
theorems 2.2 and 5.1; these replace theorems 4.2 and 5.5
from [FP1] that had been used in [FP2]. Here we perform all
steps such that all sections and homotopies are fixed on $X_0$.

The construction in proposition 5.1 [FP2] is such that
the complex $a_{*,n}$ obtained after the first $n$ steps is
already constant over the set $A^n =\cup_{j=0}^n A_j$, i.e., it
determines a \holo\ section in a \nbd\ of $A^n$, and the convergence
of these sections to $f$ is uniform on compacts in $X$ as
$n\to\infty$. The modification process also gives at each step
a homotopy (a 1-prism) between the two adjacent complexes
$a_{*,n}$ and $a_{*,n+1}$. This gives in the end
a 1-prism connecting the initial section $a=a_{*,0}$ and
the limit section $f$.

Finally, applying the version of proposition 5.1 in [FP2]
for continuous prisms, we can homotopically deform the above
1-prism from $a$ to $f$ (keepings the end sections $a$ and $f$
fixed) to another 1-prism which consists of {\it sections}
over $X$, thereby obtaining a homotopy of sections connecting
$a$ to $f$ and satisfying theorem 1.4. Complete details
of this argument can be found in sect.\ 5 and 6 of [FP2].
The same proof applies to $P$-sections.
This proves the first part of theorem 1.4.

Assume now that the initial section $a=f_0 \colon X\to Z$
is \holo\ in an open set $U_0\supset X_0\cup K$. Our goal
is to construct a \holo\ section $f\colon X\to Z$ which matches $a$
to a prescribed order $k$ along $X_0$ and is obtained from
$a$ by a homotopy $f_t\colon X\to Z$ which is \holo\ near
$X_0$ and matches $a$ to order $k$ along $X_0$.
We have two possibilities: either we replace lemma 3.2
with its analogue for functions that vanish to order $k$
along $X_0$ (see the remark following the proof of lemma 3.2),
or we apply theorem 5.2 in place of theorem 5.1 (and we don't
do any gluing along $X_0$). The first approach works exactly
as before. We now elaborate the second approach which has the
advantage that no spray is needed on $Z$ over neighborhoods
of points $x\in X_0$. The main step is the following proposition
which is similar to proposition 8.2 in [FP2].

\proclaim 6.1 Proposition:
Assume that the submersion $h\colon Z\to X$ admits a spray
over an open \nbd\ of each point $x\in X\bs (X_0\cup K)$.
Let $f_0\colon X\to Z$ be a continuous section which is
\holo\ in an open set $U_0\supset X_0\cup K$.
For each compact \hc{X}\ subset $L\subset X$ containing $K$
there are an open set $U'_0\supset X_0\cup K$ and a homotopy
of continuous sections $f_t\colon X\to Z$ $(t\in [0,1])$, with
$f_0$ being the given section, such that for each
$t\in [0,1]$, $f_t$ is \holo\ in $U'_0$, it agrees with
$f_0$ to order $k$ at $X_0$, $f_t$ approximates $f_0$
uniformly on $K$, and the section $f_1$ is \holo\ in an open
set $W_0\supset X_0\cup L$.

Granted proposition 6.1, we complete the proof of theorem 1.4
by exhausting $X$ by a sequence of compact, \hc{X}\ subsets
$K=L_0\subset L_1\subset L_2\subset \ldots$ and
applying proposition 6.1 inductively.
\endpr

\demo Proof of proposition 6.1:
Let $(A_0,A_1,\ldots,A_n)$ be a Cartan string in $X$
for the pair of sets $K\subset L$, provided by lemma 8.4
in [FP2], satisfying

\item{(i)}   $K\cup (X_0\cap L) \subset A_0 \subset U_0$;

\item{(ii)}  for $j=1,2,\ldots,n$ we have
$A_j \cap (X_0 \cup K)=\emptyset$ and $A_j \subset U_k$
for some $k=k(j)\ge 1$;

\item{(iii)} $L=\cup_{0\le j\le n} A_j$.
\medskip\rm

By (ii) we have $A_j\cap X_0=\emptyset$ for all $j\ge 1$
and hence this string is $X_0$-regular. We now
inductively apply proposition 5.1 in [FP2], making
certain to satisfy the interpolation condition along $X_0$.
In the process we must patch \holo\ sections on Cartan pairs
$(A^j,A_{j+1})$, where $A^j=A_0\cup A_2\cup\ldots \cup A_j$,
as well as on some other Cartan pairs which do not intersect $X_0$.
While no change in the procedure is needed for the second
case, the first case must be modified in such a way that
the resulting sections will be \holo\ in a \nbd\ of $X_0$
and will interpolate the initial section.
Set $A'_0= A_0\cup X_0$ and $\cA'=(A'_0,A_1,\ldots,A_n)$.
For each $j$ the set $A^j=\cup_{l=0}^j A_l$ is \hc{X}, and
hence the same is true for
$(A')^j = X_0\cup A^j = A'_0\cup A_1\cup\ldots \cup A_j$.
Thus $\cA'$ satisfies all properties of an $X_0$-regular
Cartan string, except that the initial set $A'_0$ is unbounded.

By proposition 4.7 in [FP2] we can associate to $f_0$ the initial
holomorphic $\cK(\cA')$-complex $f_{*,0}$ such that the section
associated to $A'_0$ equals the restriction of $f_0$
to a suitable open \nbd\ of $A'_0$ in $X$. We now modify this
complex $f_{*,0}$ by inductively applying proposition 5.1 in [FP2],
except that we perform the gluing of sections on Cartan pairs
$((A')^j,A_{j+1})$ by appealing to theorem 5.2 instead of
theorem 5.1. The \nbd\ of $A'_0$ on which the new section
(the result of the gluing) is defined may shrink at each
step without any detrimental consequences. After $n$ steps
we obtain a \holo\ section $f_1$ in an open set
$W'_0 \supset X_0\cup L$ which can be connected
to $f_0$ by a \holo\ homotopy $f_t\colon W'_0 \to Z$
$(t\in [0,1])$.

It remains to extend $f_t$ to $X$ by modifying it outside
a smaller \nbd\ of $X_0\cup L$. This is done as usual by taking
$f_{t\chi(x)}(x)$ $(x\in X,\ t\in [0,1])$,
where $\chi\colon X\to [0,1]$ is a smooth cut-off
function which vanishes outside $W'_0$ and equals one
in a smaller open set $W_0\supset X_0\cup L$.
This completes the proof.
\endpr

%
%   Stratified sprays
%
\beginsection 7. Submersions with stratified sprays.

We can extend theorem 1.4 to submersions with stratified
sprays over Stein spaces (compare with sec.\ 3.1 in
Gromov [Gro]), a situation which arises naturally
in global analytic geometry. An important special case
was considered in 1966 by Forster and Ramspott [FR].

Let $X$ and $Z$ be complex spaces, possibly
with singularities. A \holo\ map $h\colon Z\to X$
is said to be a submersion (of co-rank $k$) if it is
locally near each point $z_0\in Z$ equivalent
(by a fiber preserving biholomorphic map) to a projection
$p\colon U\times V\to U$, where $U\subset X$ is an open set
in $X$ containing $x_0=h(z_0)$ and $V$ is an open set in
some $\C^d$.

Assume now that the base $X$ of a submersion $h\colon Z\to X$
is a Stein space which is stratified by a descending chain
of closed complex subspaces
$X=X_m\supset X_{m-1}\supset\cdots\supset X_1\supset X_0$
such that each stratum $Y_k=X_k\bs X_{k-1}$ ($1\le k\le m$)
is regular and the restricted submersion
$h\colon Z|_{Y_k} \to Y_k$ admits a spray over a small
\nbd\ of any point $x\in Y_k$ (def.\ 2 in sec.\ 1).
Then the following holds:

\medskip
\ni {\it Any continuous section $f_0\colon X\to Z$ such that
$f_0|_{X_0}$ is \holo\ can be deformed to a \holo\
section $f_1\colon X\to Z$ by a homotopy
$f_t\colon X\to Z$ $(t\in [0,1])$ that is fixed on $X_0$.}
\medskip

A similar result holds for $P$-sections whose restriction
to $X_0$ is \holo\ (see def.\ 3 and theorem 1.4).

One can prove this by induction over the strata $X_k$ as follows.
Suppose that for some $k\ge 1$ we have already constructed a
continuous section $f_{k-1}\colon X\to Z$ that is \holo\ over $X_{k-1}$
and satisfies $f_{k-1}=f_0$ on $X_0$. We first apply proposition 1.3
(which holds also for complex spaces) to make $f_{k-1}$
\holo\ in an open \nbd\ of $X_{k-1}$ in $X$ and continuous
elsewhere. We then apply theorem 1.4, with $(X,X_0)$ replaced by
the pair $(X_k,X_{k-1})$, to deform $f_{k-1}$ through
a homotopy to a continuous section $f_k\colon X\to Z$
that is \holo\ over $X_k$ and matches $f_0$ on $X_0$.
In finitely many steps we obtain a \holo\ section over $X$.

The construction $f_k$ requires a version of theorem 1.4
in which $X$ is a Stein space whose singular locus is contained
in the subvariety $X_0$ where the initial section is already \holo.
Such a result can be proved by small changes in the proof
of proposition 4.2 and theorem 5.2 as follows.
Each Stein space $X$ of dimension $n$ admits a holomorphic
map $g\colon X\to \C^{2n+1}$ which is a homeomorphism of $X$
onto a closed Stein subspace $\tilde X =g(X)\subset \C^{2n+1}$
and is an injective immersion on the regular part of $X$ [GuR].
Observe that the attaching of pairs of \holo\ sections
(proposition 4.2 and theorem 5.2) needs to be carried out
only on Cartan pairs $(A,B)$ in $X$ for which the set $C=A\cap B$
is contained in the regular part of $X$. Since $g$ is regular and
injective there, we may transfer the data for the patching problem
(the transition map $\psi$ in proposition 4.2 and theorem 5.2)
to a suitably chosen open \nbd\ of $g(C)$ in $\C^{2n+1}$.
The solution obtained on a Cartan pair in $\C^{2n+1}$ is then
pulled back to $X$ by $g$. The globalization scheme in sect.\ 6
goes through without changes. We leave out further details.

\demo Remark:
It seems that the first examples of
`submersion with stratified sprays'
satisfying the h-principle were the {\it Endromisb\"undel}
of Forster and Ramspott [FR]. These are submersions $Z\to X$,
where $Z$ is an open subset in $X\times M(s,r)$
($M(r,s)=$the set of complex $s\times r$ matrices)
and the restrictions $Z|_{Y_k}\to Y_k$ to certain strata
$Y_k=X_k\bs X_{k-1}$ are \holo\ fiber bundles with complex
Grassman manifolds as fibers. The corresponding h-principle,
proved in [FR], had many interesting applications, e.g.\ to
obtain the minimal number of global generators of a given
coherent analytic sheaf over $X$.

\demo Acknowledgements: We wish to thank Bo Berndtsson for his
help in the proof of lemma 3.1. The first author acknowledges
partial support by the NSF, by the Vilas foundation at the
University of Wisconsin--Madison and by the Ministry of Science
of the Republic of Slovenia. The second author was supported
in part by the Ministries of Science and of Education in Slovenia.
A part of the work was done while the second author was visiting 
the Mathematics Department at the University of Wisconsin--Madison 
in the Spring of 2000 and she wishes to thank this institution 
for its hospitality.

%
%
%  References
%
%
\medskip\ni\bf References. \rm

\ii{[AAC]} K.\ Adachi, M.\ Andersson, H.R.\ Cho:
$L\sp p$ and $H\sp p$ extensions of holomorphic functions from
subvarieties of analytic polyhedra.
Pacific J.\ Math.\ {\bf 189} (1999), 201--210.

\ii{[Car]} H.\ Cartan: Espaces fibr\'es analytiques.
Symposium Internat.\ de topologia algebraica, Mexico, 97--121 (1958).
(Also in Oeuvres, vol.\ 2, Springer, New York, 1979.)

\ii{[Dem]} J.-P.\ Demailly:
Un exemple de fibr\'e holomorphe non de Stein \`a fibre $\bf C^2$ ayant
pour base le disque ou le plan.
Invent.\ Math.\ {\bf 48}, 293--302 (1978).

\ii{[FR]} O.\ Forster and K.\ J.\ Ramspott:
Analytische Modulgarben und Endromisb\"undel.
Invent.\ Math.\ {\bf 2}, 145--170 (1966).

\ii{[FP1]} F.\ Forstneri\v c and J.\ Prezelj:
Oka's principle for holomorphic fiber bundles with sprays.
Math.\ Ann.\ {\bf 317} (2000), 117-154.

\ii{[FP2]} F.\ Forstneri\v c and J.\ Prezelj:
Oka's principle for holomorphic submersions with sprays.
Preprint, 1999.
% Math.\ Ann., to appear.

\ii{[Gr1]} H.\ Grauert:
Approximationss\"atze f\"ur holomorphe Funktionen mit Werten in
komplexen R\"aumen.
Math.\ Ann.\ {\bf 133}, 139--159 (1957).

\ii{[Gr2]} H.\ Grauert:
Holomorphe Funktionen mit Werten in komplexen Lieschen Gruppen.
Math.\ Ann.\ {\bf 133}, 450--472 (1957).

\ii{[Gr3]} H.\ Grauert: Analytische Faserungen über
holomorph-vollständigen Räumen.
Math.\ Ann.\ {\bf 135}, 263--273 (1958).

\ii{[GRe]} H.\ Grauert, R.\ Remmert: Theory of Stein Spaces.
Grundl.\ Math.\ Wiss.\ {\bf 227}, Springer, New York, 1977.

\ii{[Gro]} M.\ Gromov:
Oka's principle for holomorphic sections of elliptic bundles.
J.\ Amer.\ Math.\ Soc.\ {\bf 2}, 851-897 (1989).

\ii{[GuR]} C.\ Gunning, H.\ Rossi:
Analytic functions of several complex variables.
Prentice--Hall, Englewood Cliffs, 1965.

\ii{[Hen]} G.\ M.\ Henkin:
Continuation of bounded holomorphic functions from submanifolds
in general position in a strictly pseudoconvex domain. (Russian)
Izv.\ Akad.\ Nauk SSSR Ser.\ Mat.\ {\bf 36} (1972), 540--567.

\ii{[HL1]} G.\ M.\ Henkin, J.\ Leiterer:
Theory of functions on complex manifolds.
Akademie-Verlag, Berlin, 1984.

\ii{[HL2]} G.\ Henkin, J.\ Leiterer:
The Oka-Grauert principle without induction over the basis dimension.
Math.\ Ann.\ {\bf 311}, 71--93 (1998).

\ii{[H\"or]} L.\ H\"ormander:
An Introduction to Complex Analysis in Several Variables, 3rd ed.
North Holland, Amsterdam, 1990.

\ii{[Pre]} J.\ Prezelj: Interpolation of Embeddings
of Stein Manifolds on Discrete Sets.
Preprint, 2000.

\ii{[RRu]} J.-P.\ Rosay, W.\ Rudin:
Holomorphic maps from $\C^n$ to $\C^n$.
Trans.\ Amer.\ Math.\ Soc.\ {\bf 310}, 47--86 (1988)

\ii{[Siu]} Y.T.\ Siu: Every Stein subvariety admits a Stein neighborhood.
Invent.\ Math.\ {\bf 38}, 89--100 (1976).

%
%
%  Addresses
%
%
\bigskip\medskip
\settabs 5\columns
\+\ \ Franc Forstneri\v c            &&& Jasna Prezelj\cr
\+\ \ Department of Mathematics      &&& IMFM, University of Ljubljana  \cr
\+\ \ University of Wisconsin        &&& Jadranska 19 \cr
\+\ \ Madison, WI 53706, USA         &&& 1000 Ljubljana, Slovenia \cr
\+                                   \cr
\+\ \ \it Current address:\rm \cr
\+\ \ IMFM, University of Ljubljana  \cr
\+\ \ Jadranska 19  \cr
\+\ \ 1000 Ljubljana, Slovenia\cr

\bye

% IGNORE THE REST
\ni \it Address: \rm Institute of Mathematics, Physics and Mechanics,
University of Ljubljana, Jadranska 19, 1000 Ljubljana, Slovenia

\bye